\documentclass[12pt]{article}
\usepackage{amsmath}
\usepackage{graphicx,psfrag,epsf}
\usepackage{enumerate}
\usepackage[font=small,format=plain,labelfont=bf,textfont=sl]{caption}
\usepackage{wrapfig}
\usepackage{threeparttable} 
\usepackage{amsthm}
\usepackage{amsfonts}

\newcommand{\blind}{0}

\begin{document}

\def\spacingset#1{\renewcommand{\baselinestretch}%
{#1}\small\normalsize} \spacingset{1}

%%%%%%%%%%%%%%%%%%%%%%%%%%%%%%%%%%%%%%%%%%%%%%%%%%%%%%%%%%%%%%%%%%%%%%%%%%%%%%

\if0\blind
{
  \title{\bf Principal Nested Spheres for Time Warped Functional Data Analysis}
  \author{Xiaosun Lu and J. S. Marron\\
    Department of Statistics and Operations Research\\
    University of North Carolina Chapel Hill}
  \maketitle
} \fi

\if1\blind
{
  \bigskip
  \bigskip
  \bigskip
  \begin{center}
    {\LARGE\bf Title}
\end{center}
  \medskip
} \fi

\bigskip
\begin{abstract}
There are often two important types of variation in functional data: the horizontal (or phase) variation and the vertical (or amplitude) variation. These two types of variation have been appropriately separated and modeled through a domain warping method (or curve registration) based on the Fisher Rao metric. This paper focuses on the analysis of the horizontal variation,  captured by the domain warping functions. The square-root velocity function representation transforms the manifold of  the warping functions to a Hilbert sphere. Motivated by recent results on manifold analogs of principal component analysis, we propose to analyze the horizontal variation via a Principal Nested Spheres approach. Compared with earlier approaches, such as approximating tangent plane principal component analysis, this is seen to be the most efficient and interpretable approach to decompose the horizontal variation in some examples.
\end{abstract}

\noindent%

\spacingset{1.45}

%================================================================
\section{Introduction}

A common variability in functional data is that prominent features in the functions vary in position from one sample to another, such as the timing variation of the adolescent growth spurt in human growth curves. This positional, or phase, variation is called the {\itshape{horizontal variation}}. Another important component of variability in functional data is the amplitude variation, or the {\itshape{vertical variation}}, such as the height difference among the individuals.

There is a large literature on statistical analysis of functions, such as Kneip and Gasser (1992) \cite{KT}, Locantore et al (1999) \cite{LM}. A general overview of functional data analysis is provided by Ramsay and Silverman (2002 \cite{ramsay2002applied} and 2005 \cite{FDA}). 
Plenty of useful tools and methods are available, such as the Functional Principal Component Analysis (FPCA), with many important applications in a wide variety of
scientific fields. One open problem in functional data analysis is that, in those traditional approaches, the functional data are analyzed under the $\mathbb{L}^2$ metric, which tends to strongly focus on the vertical variation. The horizontal variation cannot be easily understood in these vertical analyses. Section \ref{SEC:FPCA} gives examples illustrating the shortcoming of conventional FPCA. 

The main purpose of this paper is to find an improved method for the {\itshape{horizontal analysis}}, i.e. the analysis of the horizontal variation. Considering the special spherical structure of the horizontal variation (see Section \ref{SEC:hSRVF} for details),  we propose to use an approach involving Principal Nested Spheres (PNS) introduced by Jung et al (2010) \cite{pns}. Comparison with several other popular approaches, such as the FPCA, suggests improved efficiency of PNS for horizontal analysis. A toy example (the left panel of Figure \ref{FIG:align}) is used to illustrate the advantages of the PNS approach.  

We find Object Oriented Data Analysis, introduced by Wang and Marron (2007) \cite{ooda}, very helpful in studying horizontal variation. The data objects are understood as the {\itshape{atoms}} of the analysis. In this study, they are functions. Several different types of data objects (or functions) are considered in this paper for horizontal analysis. These different choices of data objects and the corresponding horizontal analyses are discussed in Section \ref{SEC:ha}.

%================================================================
%================================================================
\subsection{Function Alignment Based on the Fisher Rao Metric}

\begin{figure}[htpb]
\center
\includegraphics[keepaspectratio=true, width=0.25\textwidth,
page=1]{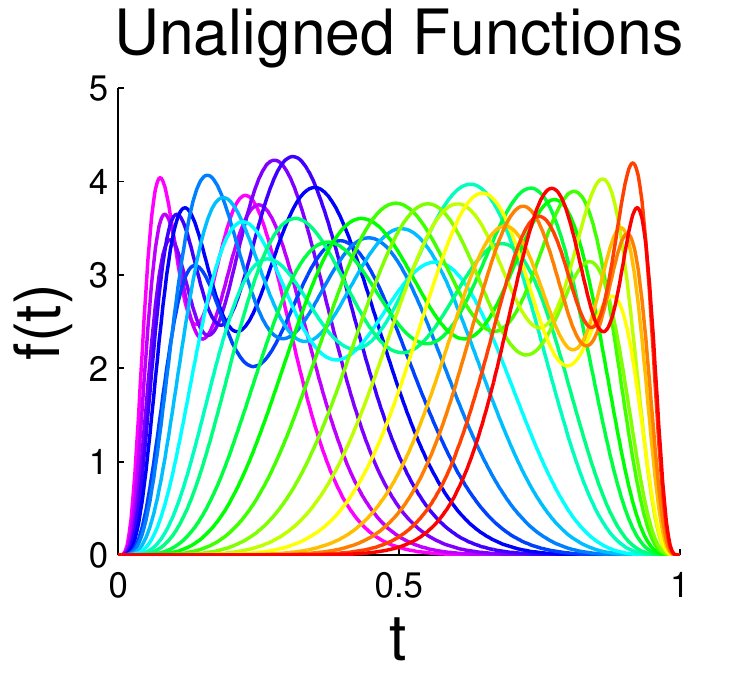} 
\includegraphics[keepaspectratio=true, width=0.25\textwidth,
page=1]{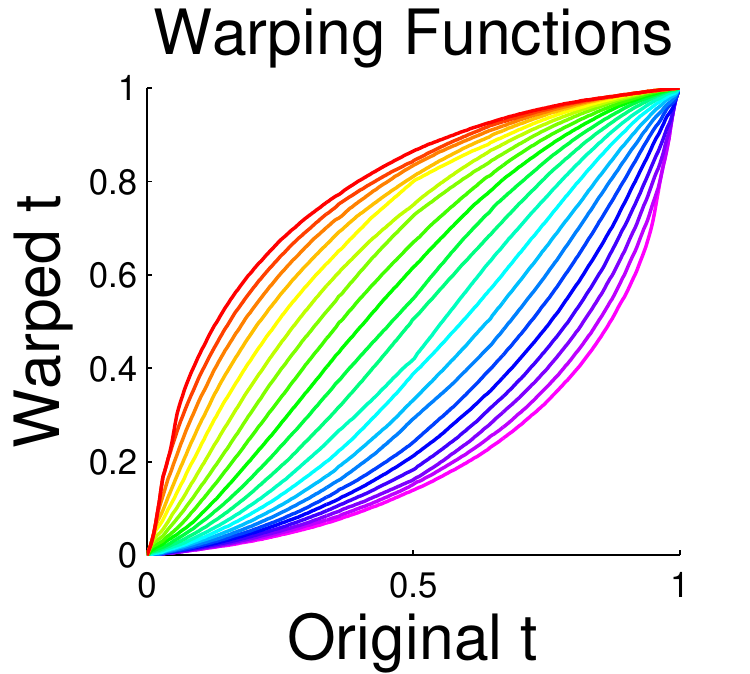}
\includegraphics[keepaspectratio=true, width=0.25\textwidth,
page=1]{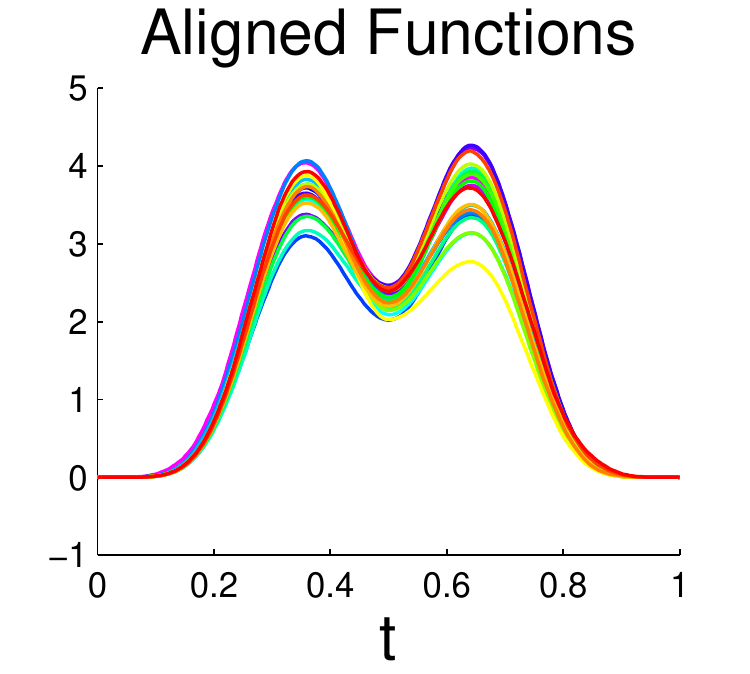}
\caption{Left: A toy example of bimodal functions with big horizontal variation. The color reflects the order of the horizontal positions of the peaks. Middle: The domain warping functions to align the functions based on the Fisher Rao metric. Right: The aligned functions.}
\label{FIG:align}
\end{figure}

A useful approach to horizontal analysis is through the idea of {\itshape{elastic functions}}. 
Some pioneering work in this area includes Ramsay and Li (1998) \cite{RL}, Gervini and Gasser (2004) \cite{GG}, Liu and Mueller (2004) \cite{FR_liu}, Kneip and Ramsay (2008) \cite{KR}, Tong and Mueller (2008) \cite{TM}.
The basic idea is to first separate the vertical and the horizontal variation through {\itshape{function alignment}}, or curve registration. In particular, consider a collection of functions $f_i(t)$, $i=1, 2, ...,n$, having both vertical and horizontal variation, such as the bimodal functions shown in Figure \ref{FIG:align} (left).
If these functions are well aligned by warping the domain properly, then the horizontal and the vertical variation can be separately captured by the domain warping functions $\gamma_i(t)$ and the resulting aligned functions $f_i(\gamma_i(t))$, respectively. For this toy example, such a set of warping functions and the corresponding aligned functions are shown in the middle and right panels respectively (details about finding those warping functions are discussed later). Then, the horizontal analysis can be done by studying those warping functions.

\begin{figure}[htpb]
\center
\includegraphics[keepaspectratio=true, width=0.6\textwidth,
page=1]{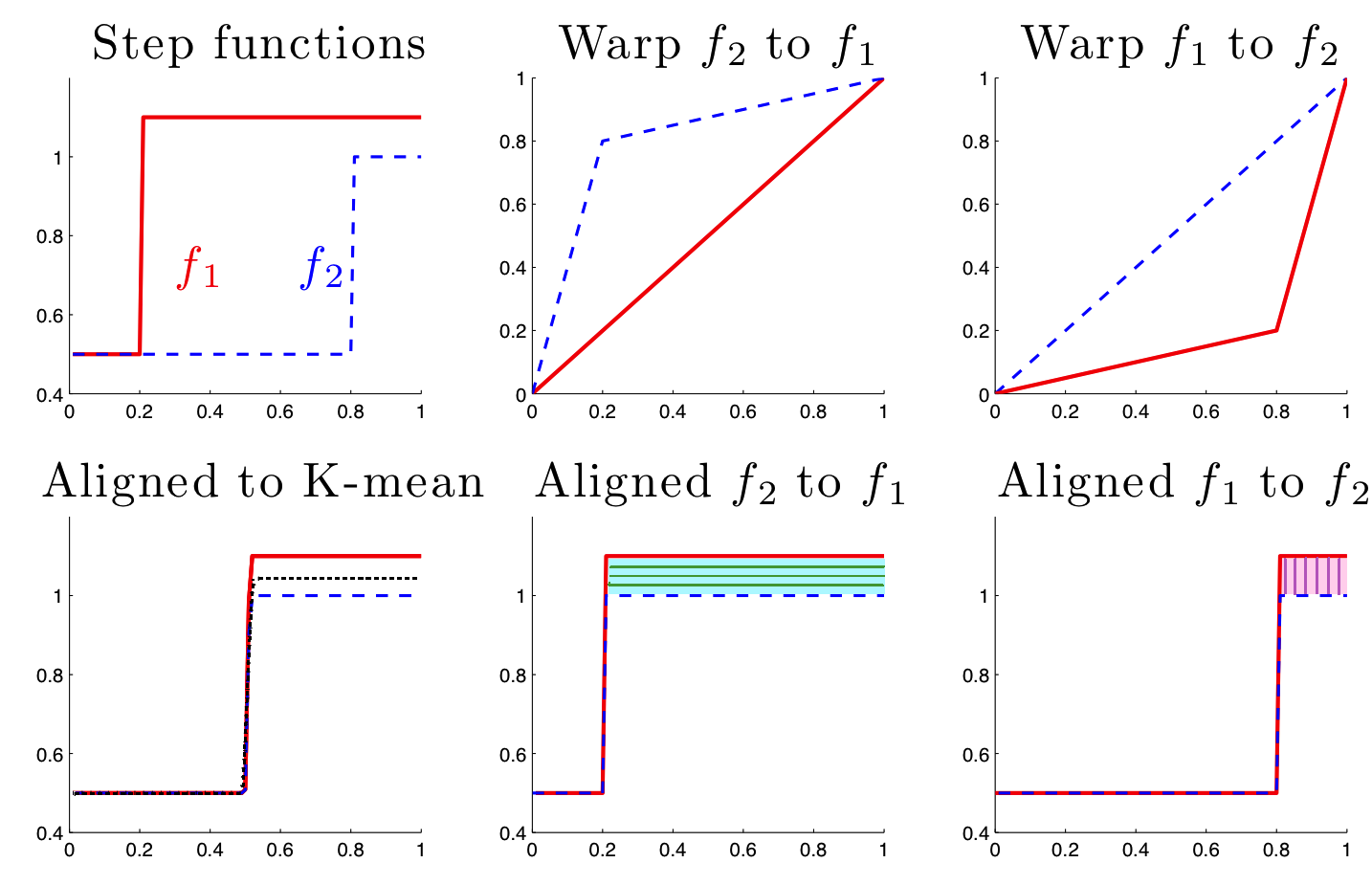}
\caption{The problem with $\mathbb{L}^2$ metric alignment. The top left panel shows two step functions $f_{1}$ (solid red) and $f_{2}$ (dashed blue). The four right panels show that warping $f_{2}$ to $f_{1}$ (middle) is different from warping $f_{1}$ to $f_{2}$ (right) under the $\mathbb{L}^2$ metric. The top two panels show the warping functions, while the bottom two panels show the aligned functions. Better than either is the Fisher Rao alignment shown in the bottom left panel, where the black dotted line indicates the Karcher mean function. }
\label{FIG:alignIssue}
\end{figure}

A crucial step in the function alignment is to find appropriate domain warping functions. Consider two functions $f_1$ and $f_2$. Most of the past approaches involve solving 
$\inf_{\gamma\in\Gamma}\|f_{1}-(f_{2}\circ\gamma)\|$
 to align $f_{2}$ to $f_{1}$, where $\|\cdot\|$ is the standard $\mathbb{L}^2$ metric, i.e.
$\|f\| = (\int_{0}^{1}|f(t)|^2 dt)^{1/2}$. However, this criterion is problematic, since the objective function  is not symmetric
in the sense that aligning $f_{1}$ to $f_{2}$ leads to a different optimal minimum. To illustrate this point, Figure \ref{FIG:alignIssue} shows a simple example of aligning two step functions. It is seen that aligning $f_{2}$ to $f_{1}$ (middle) and aligning $f_{1}$ to $f_{2}$ (right) are different under the $\mathbb{L}^2$ metric. The difference between the horizontally hatched blue area in Panel (2, 2) and the vertically hatched pink area in Panel (2, 3) indicates that the two corresponding objective functions $\|f_{1}-(f_{2}\circ\gamma)\|$ and $\|(f_{1}\circ\gamma)-f_{2}\|$ are not equal. This is because the $\mathbb{L}^2$ metric is not invariant under re-parameterization, or domain warping. In particular, $\|f_{1}-f_{2}\|\neq \|f_{1}\circ\gamma - f_{2}\circ\gamma\|$. A more appropriate metric is the Fisher Rao metric. See Srivastava et al (2011) \cite{FR} for definition and relevant theory. This metric is derived from a Riemannian metric first introduced by C. R. Rao (1945) \cite{FRdist}. A nice property of the Fisher Rao metric is that it is warping-invariant. In fact, Cencov (1982) \cite{FRdist4} proved that it is the only metric that has this property. Thus, we propose to use the Fisher Rao metric to align functions for the purpose of the horizontal analysis.

Calculations based on the Fisher Rao metric are difficult. In practice, a convenient square-root velocity function (SRVF) representation, i.e. transforming the function $f(t)$ to $\frac{\dot{f}(t)}{\sqrt{\mid\dot{f}(t)\mid}}$, simplifies the Fisher Rao framework. Under the SRVF representation, the Fisher Rao metric becomes the standard $\mathbb{L}^2$ metric, and thus, standard statistical tools for the $\mathbb{L}^2$ space, such as mean, covariance and principal components, can be used.
As an example of function alignment based on the Fisher Rao metric, the warping functions (middle) for the toy data (left) in Figure \ref{FIG:align} are found by an automatic and unsupervised approach based on this metric, proposed by Srivastava et al (2011) \cite{FR}.

%================================================================
\subsection{PNS for Spherical Structure of Horizontal SRVFs}\label{SEC:hSRVF}
One major benefit of the SRVF representation is that it transforms the manifold of the warping functions to a Hilbert sphere. In particular, suppose the domain of the functions $f_i$, $i=1, 2, ...,n$, to be $[0,1]$ (if not, consider a linear transformation that maps the domain to $[0,1]$). Let $\gamma_i$ be a warping function for $f_i$, i.e. $\gamma_i \in\Gamma=\{ \gamma:\ [0,1]\rightarrow[0,1]|\ \gamma(0)=0,\ \gamma(1) = 1$ and $\gamma$ is a diffeomorphism\}, where a {\itshape{diffeomorphism}} refers to a bijective differentiable function whose inverse is also differentiable. 
Then the SRVF $\psi_i$ of the warping function $\gamma_i$, referred to later as the {\itshape{horizontal SRVF}}, can be written as $\sqrt{\dot{\gamma_{i}}}$. Noting that $\|\psi_i\|^2 = \int_{0}^{1}\psi_i(t)^2 dt = \int_{0}^{1}\dot{\gamma_i}(t) dt = \gamma_i(1)-\gamma_i(0)=1$, these horizontal SRVFs naturally lie on the surface of a Hilbert unit sphere.

Due to the spherical structure of the horizontal SRVFs, we propose to use the PNS method for horizontal analysis. Introduced by Jung et al (2010) \cite{pns}, PNS is an extension of PCA for curved manifolds, especially for high dimensional spheres. This method finds a sequence of subspheres that best approximates the data using a backward approach, which starts with the high dimensional sphere and finds the best fitting subsphere of one dimension lower at each step. See Marron et al (2010) \cite{BackwardPCA} for more discussion of backward PCA. It has been shown in a number of cases that PNS can provide more effective analysis of manifold data than many other analogous approaches. See Pizer et al (2011) \cite{PGAskeletal} for such an example in the study of 3D shapes. 

For comparison purposes, another popular approach for data lying in curved manifolds, Principal Geodesic Analysis (PGA) proposed by Fletcher et al (2004) \cite{PGA}, is also investigated in this paper. 
Unlike PNS, it is a forward approach, starting with the {\itshape{Karcher mean}} (also called fr\'{e}chet or geodesic mean). 
In the following horizontal analysis, the Karcher mean refers to the representer defined in Definition 3 of Srivastava et al (2011) \cite{FR}.
PGA approximates the spherical surface by a tangent hyperplane centered at the Karcher mean. By performing PCA on this tangent plane, PGA finds the principal geodesics (i.e. great spheres) passing through the mean that best fit the data. 

In contrast to PGA, the PNS method finds the best fitting subsphere regardless of whether it is a great sphere or not. When the major variance is non-geodesic, PNS tends to find the best-fitting small spheres instead of only great spheres. Thus, when the data variability on the sphere is big enough, the PNS can give a much more effective decomposition of this variability than PGA. On the other hand, if the data variability is small, the PNS method does not improve much over the PGA method. This is because in this case the data do not have much curvature and can be approximated by a tangent plane well enough. In the following discussion, we focus on examples with big horizontal variation.

%================================================================
\section{Horizontal Analyses}\label{SEC:ha}

This section compares different horizontal analyses with the toy example in Figure \ref{FIG:align} (left), where the functions have big horizontal variation and the PNS method gives very useful improvement over PGA.

\begin{figure}[htpb]
\center
\includegraphics[keepaspectratio=true, width=0.25\textwidth,
page=1]{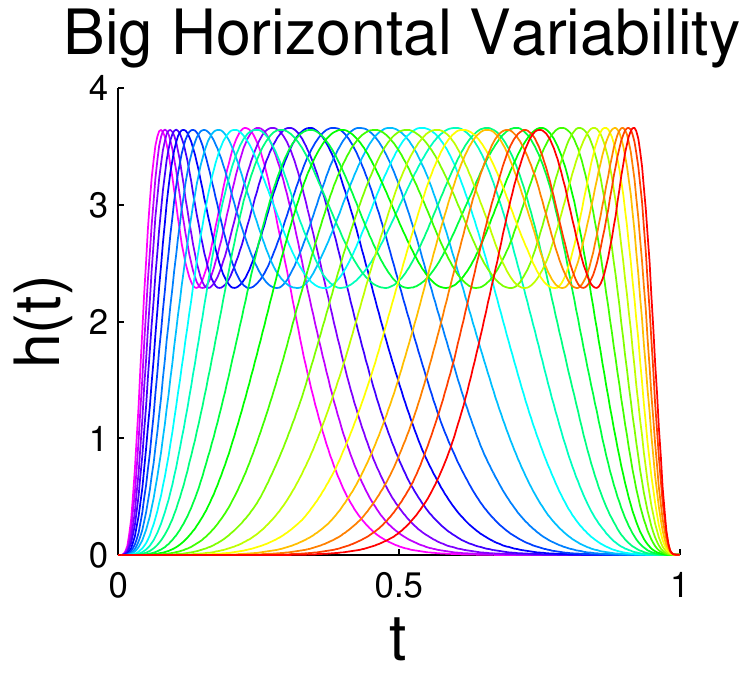} 
\includegraphics[keepaspectratio=true, width=0.25\textwidth,
page=1]{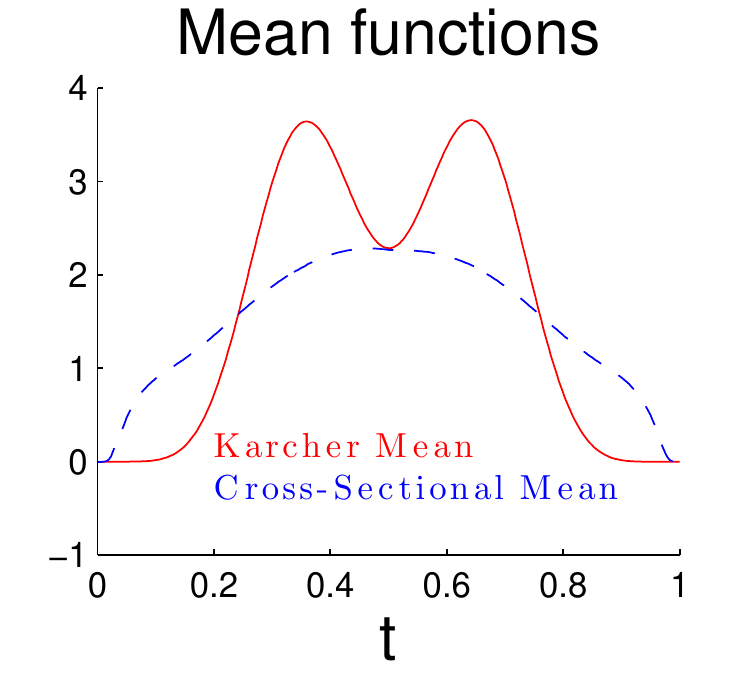}
\caption{Left: Horizontal variation of the toy example in Figure \ref{FIG:align} (left). The color reflects the order of the horizontal positions of the peaks. Right: The Karcher mean function (red solid line) and the cross-sectional mean (blue dashed line).}
\label{FIG:toyExample}
\end{figure}

Figure \ref{FIG:toyExample} (left) visualizes the pure horizontal shifts of the peaks for this toy example, which is done via warping the Karcher mean function (red curve in the right panel) by the Fisher Rao warping functions in Figure \ref{FIG:align} (middle). These functions will be referred to later as the {\itshape{horizontally shifted functions}}, denoted by $h_i$, $i=1, 2, ...,n$. 

Three different types of data objects for horizontal analysis are studied in this section: the horizontally shifted functions $h_i$, the warping functions $\gamma_i$ and the horizontal SRVFs $\psi_i$. Oriented by the choice of data objects, four different horizontal analyses have been performed on the toy data: (1) FPCA of the horizontally shifted functions $h_i$; (2) FPCA of the warping functions $\gamma_i$; (3) PGA of the horizontal SRVFs $\psi_i$; (4) PNS of the horizontal SRVFs $\psi_i$. The first two approaches, using the conventional FPCA, are discussed in Section \ref{SEC:FPCA}. 
The latter two manifold approaches, motivated by the spherical structure of the horizontal SRVFs, are discussed in Section \ref{SEC:SRVF}. Section \ref{SEC:conclusion} summarizes the comparison of these four approaches.

%================================================================
\subsection{Conventional FPCA}\label{SEC:FPCA}
An intuitive way to understand the horizontal variation of the toy data is to analyze either the horizontally shifted functions $h_i$ or the warping functions $\gamma_i$. As the FPCA is one of the most widely used statistical tools for functional data analysis, this section discusses FPCA of the two different types of functions respectively. It is seen that the conventional FPCA is rarely a good option for horizontal analysis. 

\subsubsection{FPCA of Horizontally Shifted Functions} 

\begin{figure}[htbp]
 \center
\includegraphics[keepaspectratio=true,width=0.24\textwidth,page=1]{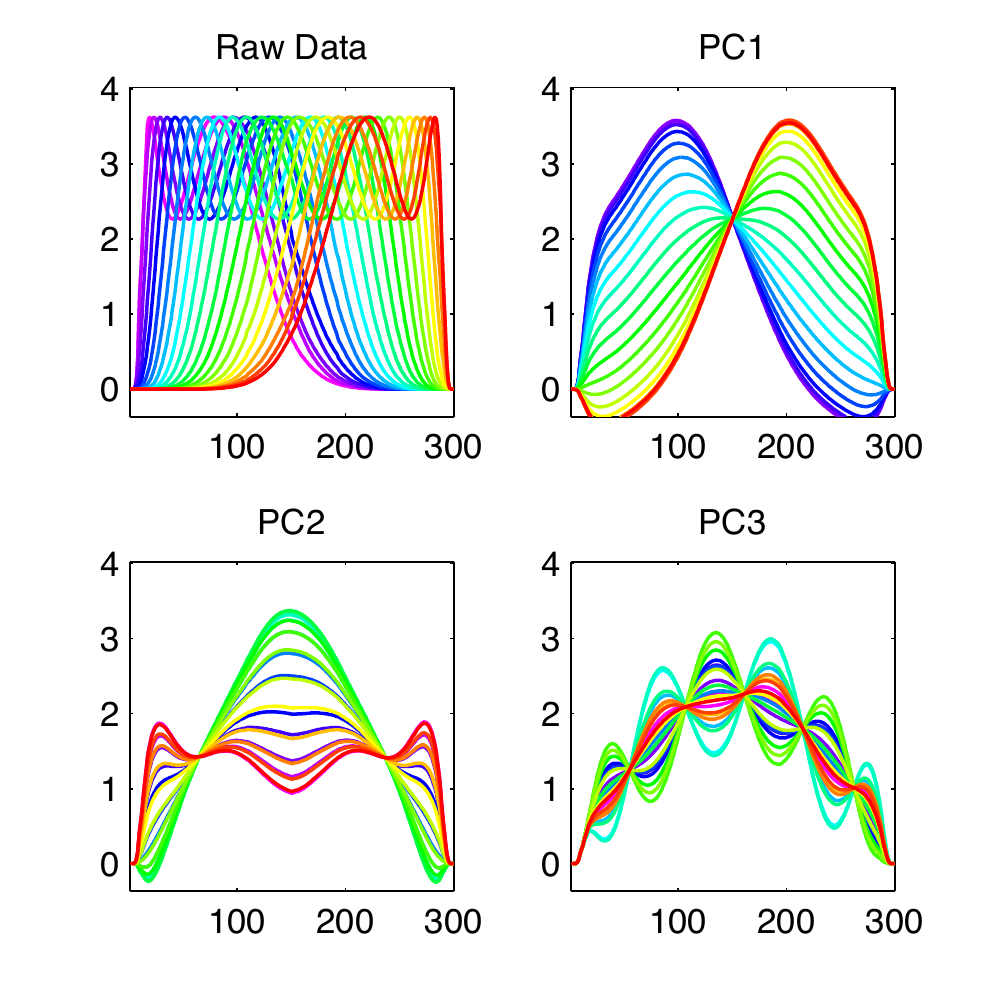}
\includegraphics[keepaspectratio=true,width=0.24\textwidth,page=1]{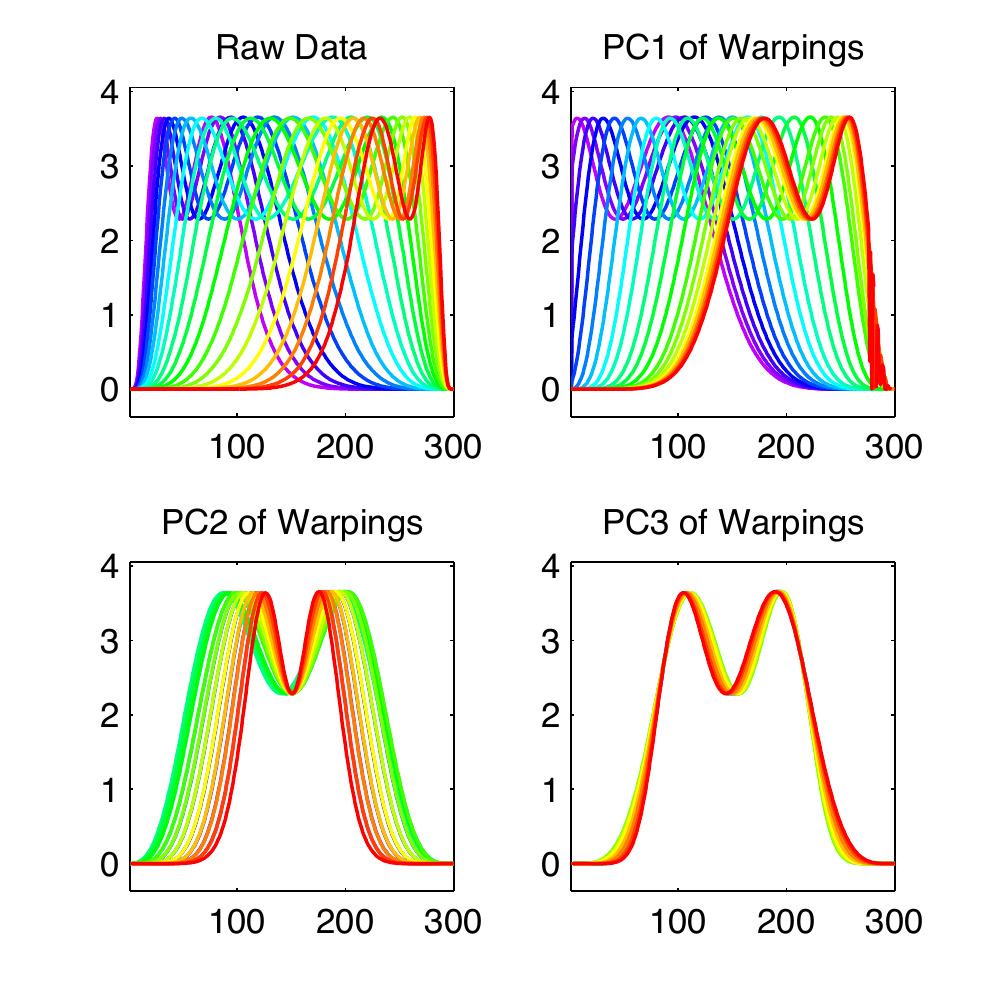}
\includegraphics[keepaspectratio=true,width=0.24\textwidth,page=1]{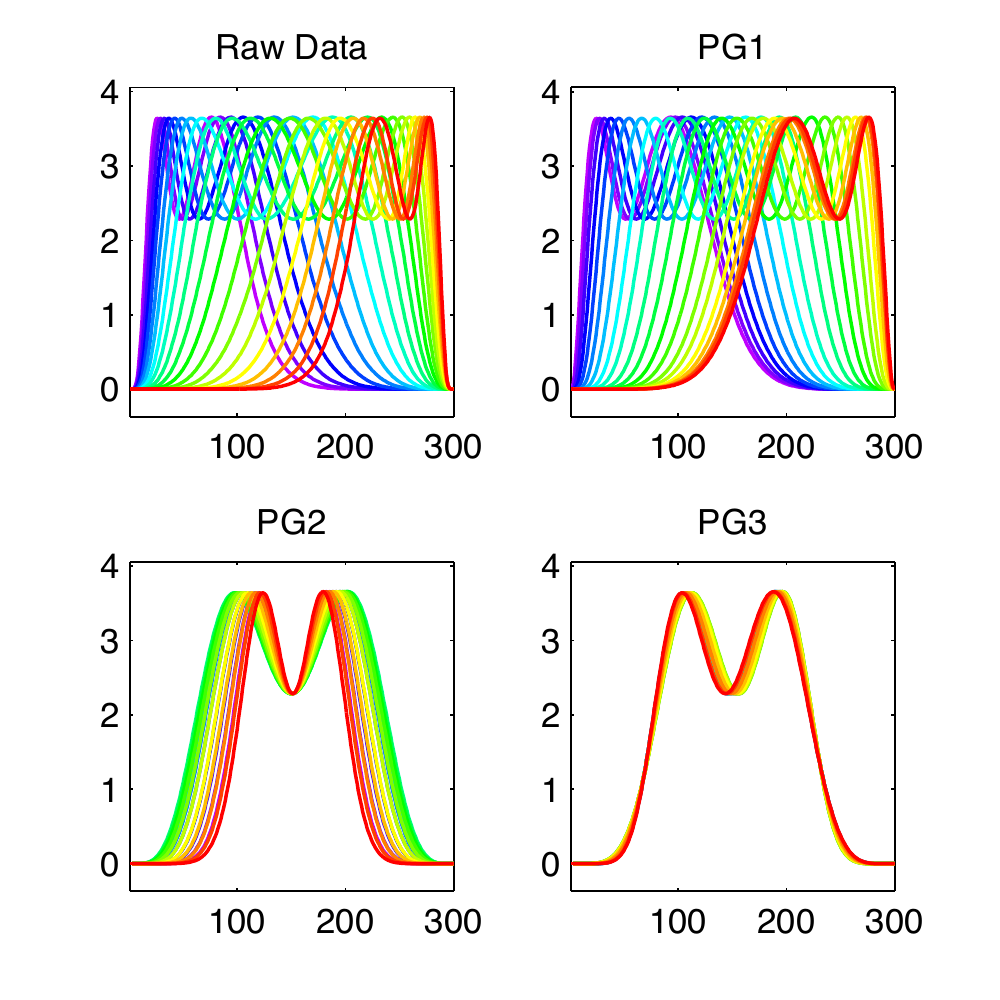}
\includegraphics[keepaspectratio=true,width=0.24\textwidth,page=1]{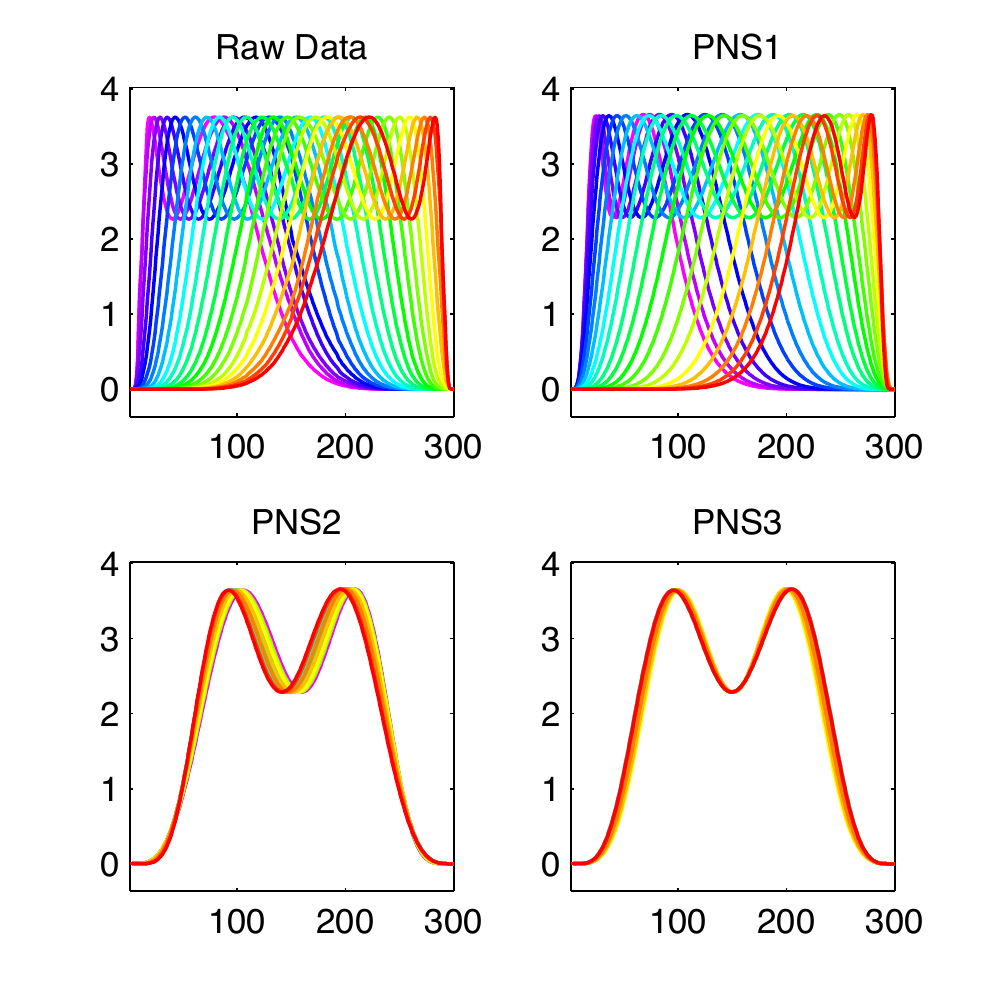}\\
\includegraphics[keepaspectratio=true,width=0.24\textwidth,page=1]{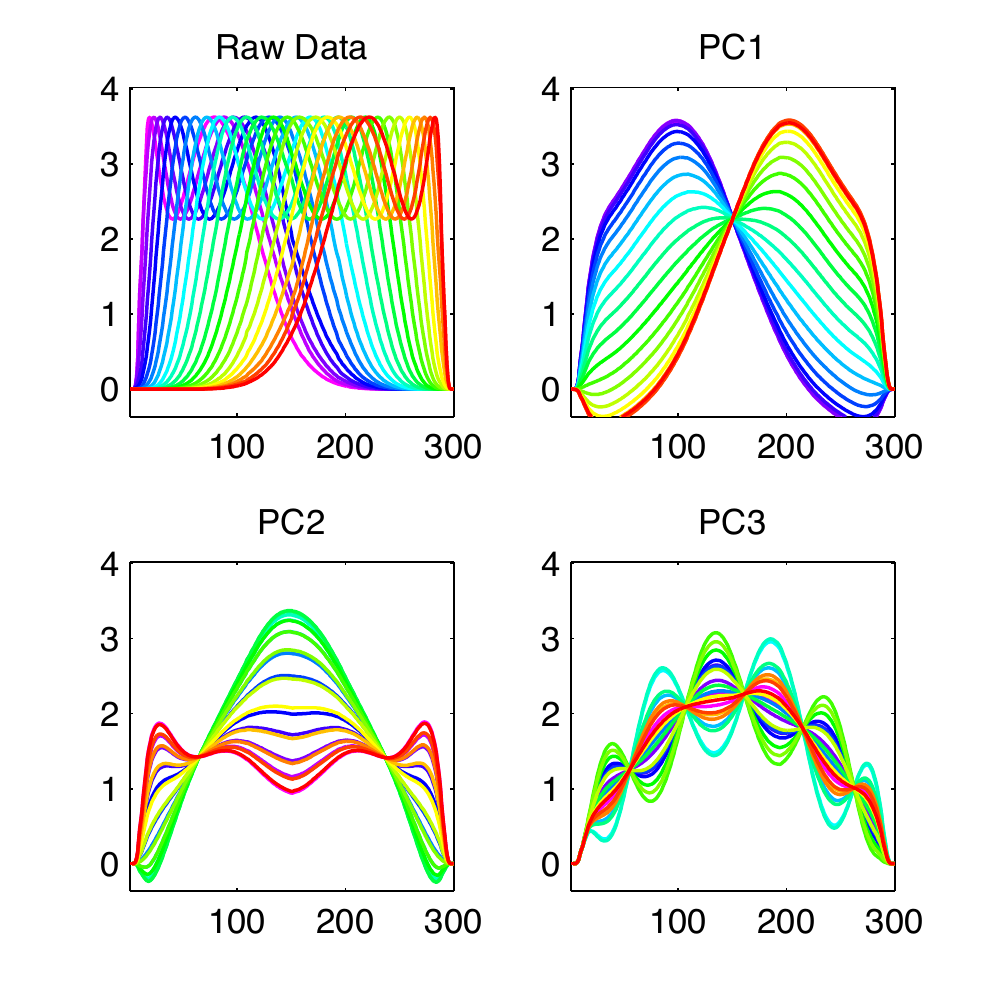}
\includegraphics[keepaspectratio=true,width=0.24\textwidth,page=1]{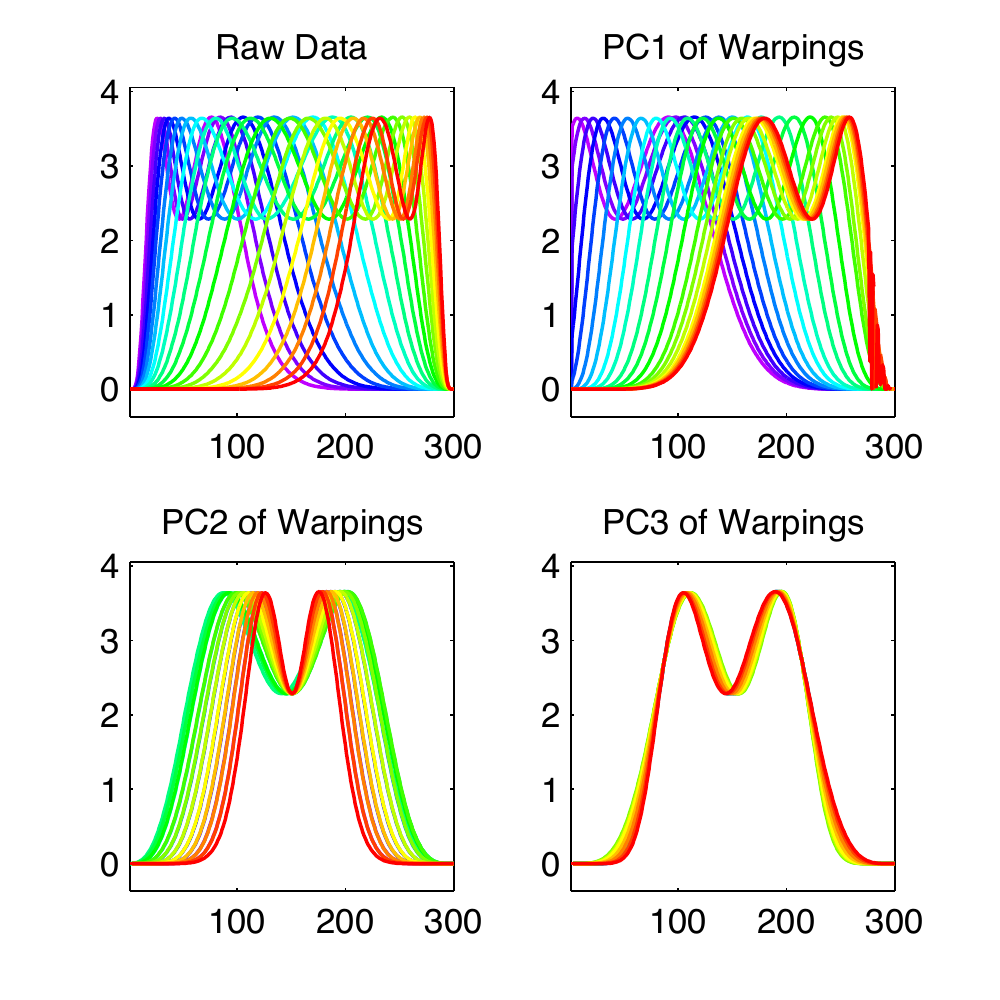}
\includegraphics[keepaspectratio=true,width=0.24\textwidth,page=1]{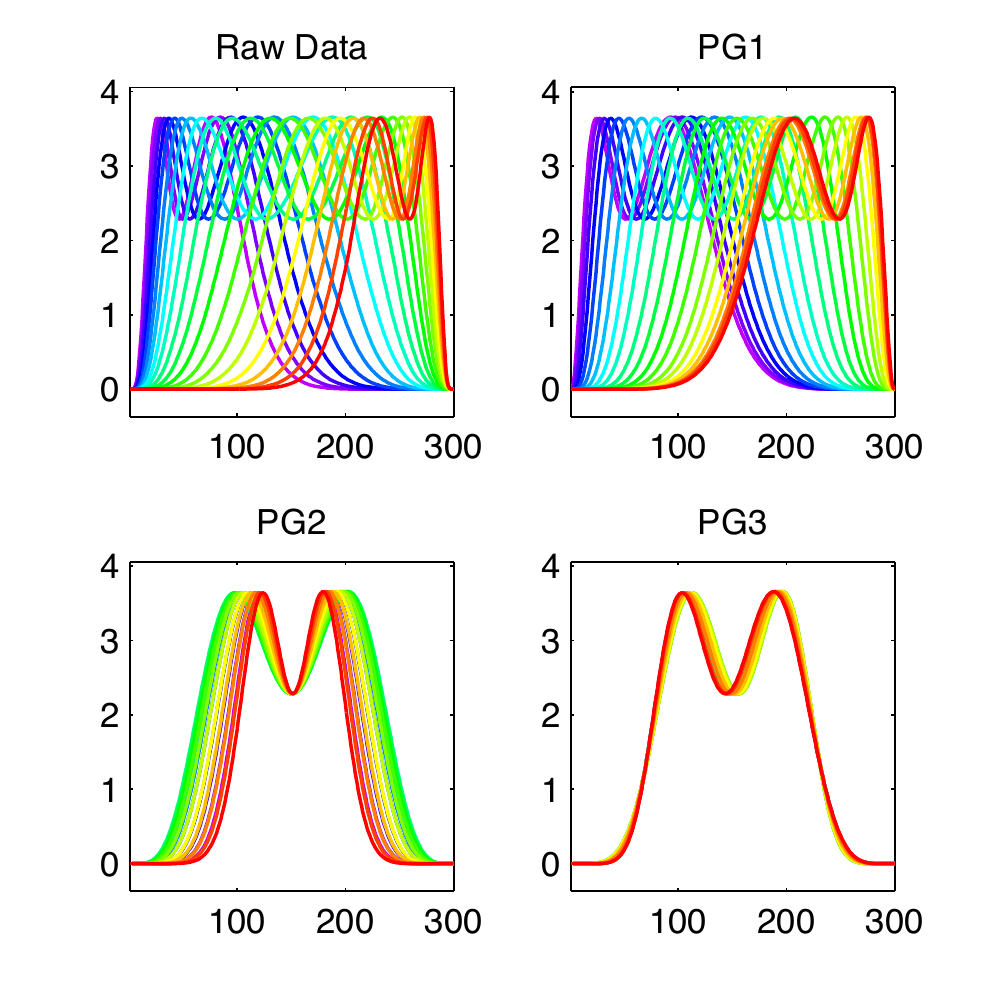}
\includegraphics[keepaspectratio=true,width=0.24\textwidth,page=1]{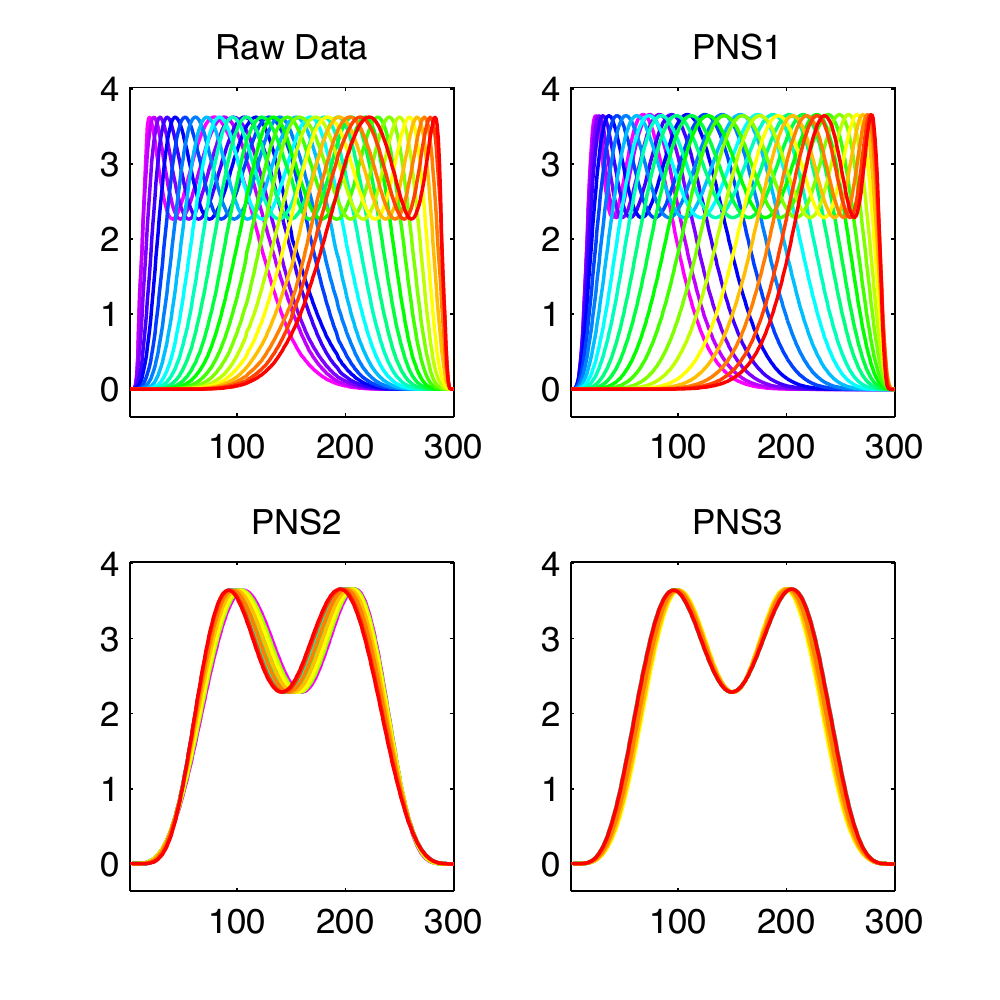}
\caption{Horizontal analyses of the toy data. The color is consistent with that in the first panel in Figure \ref{FIG:toyExample}. From the left column to the right:  FPCA of the horizontally shifted functions, FPCA of the warping functions, PGA of the horizontal SRVFs, PNS of the horizontal SRVFs. 
Each column shows the first two components of each analysis. Note that successive improvement in quality of data representation and signal compression are shown.}
\label{FIG:proj}
\end{figure}

FPCA involves centering data with the cross-sectional mean based on the $\mathbb{L}^2$ metric. However, in the case of big horizontal variability, this cross-sectional mean may hardly be useful in capturing the underlying shape of the functions. In the toy example, the unimodal cross-sectional mean (the blue dashed line in the right panel of Figure \ref{FIG:toyExample}) is a poor representative of the bimodal functions $h_i$, while the bimodal Karcher mean (the red solid line) based on the Fisher Rao metric is much more satisfying. Thus, the resulting PC projections from the FPCA of the horizontally shifted functions $h_i$ are difficult to interpret. See panels in the first column of Figure \ref{FIG:proj} for the first two PC projections. Therefore, the FPCA of the $h_i$ is not an appropriate approach for horizontal analysis.  

%=====================================
\subsubsection{FPCA of Warping Functions}\label{SEC:FPCAWarp}

One challenge of performing the FPCA on the domain warping functions is how to interpret each component. Better interpretation comes from transforming the decomposition of the warping functions into the original function space, i.e. warping the Karcher mean function by the PC projections.  
The second column of Figure \ref{FIG:proj} shows the first two transformed PC projections for the toy example. These two components provide a much more useful summary of the apparent horizontal variation in the raw data than the previous ones from the FPCA of the $h_i$ (first column).  The first component reflects the horizontal shifts of the peaks, while the second one is about the horizontal distance between the two peaks.

\begin{figure}[htbp]
\center
\includegraphics[keepaspectratio=true, scale=0.35,page=1]{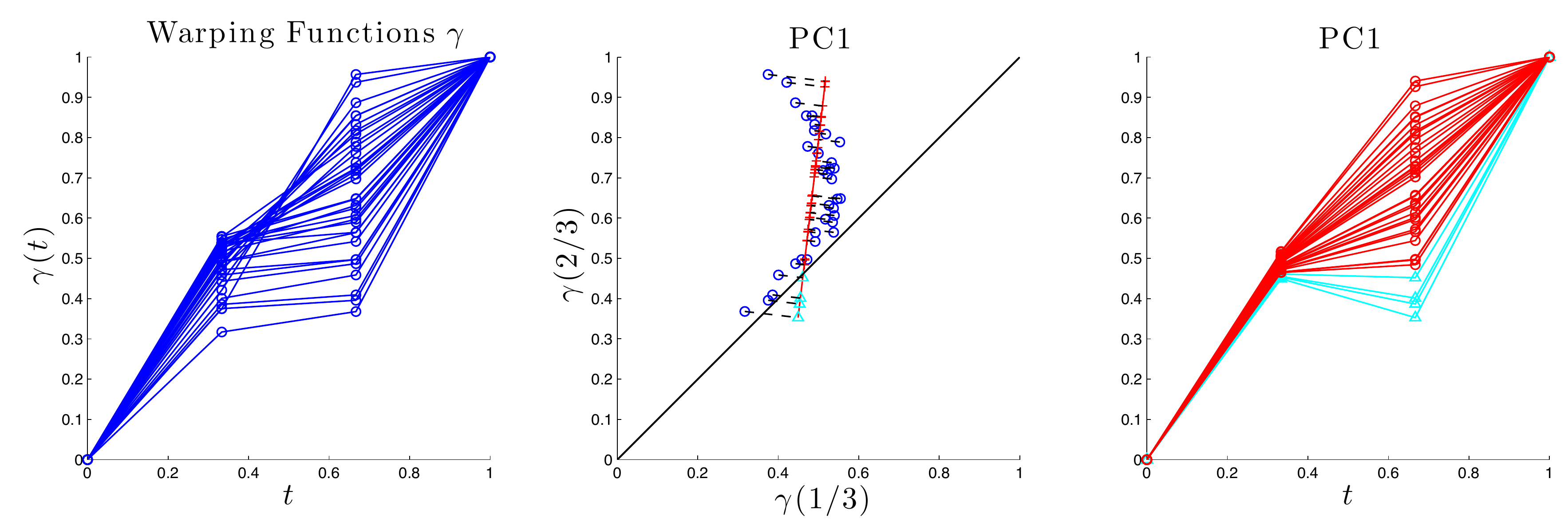} 
\caption{A toy example to illustrate the problem of performing the FPCA on warping functions. Left: A set of two-dimensional warping functions $\gamma_i$, each determined by $\gamma_i(1/3)$ and $\gamma_i(2/3)$. Middle: The scatter plot of $\gamma_i(1/3)$ and $\gamma_i(2/3)$ (blue circles), and the PC1 direction (red line) of these points. The red crosses indicate the PC1 projections above the black diagonal line, and the cyan triangles indicate the PC1 projections below the line.  Right: The projected curve visualization of those PC1 projections. Note that the cyan curves are not bijective, i.e. not valid warping functions.}
\label{FIG:warpIssue}
\end{figure}

However, this approach has a serious weakness. That is, the PC projection of a warping function is not necessarily bijective, and thus, not a warping function. In other words, the conventional FPCA leaves the space $\Gamma$ of warping functions. 
To illustrate this, Figure \ref{FIG:warpIssue} shows the FPCA of a set of simple two-dimensional warping functions $\gamma_i$ (left panel), each of which is determined by two values, $\gamma_i(1/3)$ and $\gamma_i(2/3)$. It is seen in the right two panels that some of the PC1 projections (cyan) have a decreasing part, i.e. $\gamma_i(1/3)>\gamma_i(2/3)$. 
Warping the Karcher mean function with these non-warping PC projections is problematic. Computationally, this causes the wiggly right end of the yellow and the red functions in Panel (1, 2) of Figure \ref{FIG:proj}.
 
%================================================================

\subsection{Analyses on SRVF Manifold}\label{SEC:SRVF}
The following analyses avoid the problem shown in Section \ref{SEC:FPCAWarp}, by appropriately using the spherical structure of the horizontal SRVFs $\psi_i$. The idea is to first decompose the variability of the spherical SRVFs and then transform the projections of the SRVF components back to the warping function space $\Gamma$ using the formula 
$\gamma_{\xi}(t) = \int_{0}^{t}\xi(t)^2 dt$, where $\xi$ is a point on the SRVF sphere. It can be easily checked that $\gamma_\xi\in\Gamma$. Finally, the decomposition of the horizontal variation of the original functions can be obtained via warping the Karcher mean function with the transformed SRVF projections. The following discussion shows how manifold approaches, especially the PNS approach,  can work better for the horizontal analysis than the conventional FPCA.

%================================================================
\subsubsection{PGA of Horizontal SRVFs}

The third column in Figure \ref{FIG:proj} shows the first two components of the horizontal variation in the toy data, based on the PGA of the horizontal SRVFs. Compared with the previous FPCA results (the first two columns), this approach gives a much better decomposition of the horizontal variation. 
The first component captures most of the horizontal shifts of the peaks. The second one is similar to the PC2 of the warping functions in Panel (2, 2), but has a visually smaller horizontal variation. This shows more of the underlying signal in the data has been moved to PG1, i.e. better {\itshape{signal compression}}.

%================================================================
\subsubsection{PNS of Horizontal SRVFs}
The first two components of the horizontal variation in the toy data based on the PNS of the horizontal SRVFs are shown in the fourth column in Figure \ref{FIG:proj}. Results from this decomposition give more signal compression than those from the previous analyses. 
The first component simultaneously captures both the mode of peak location and the mode of distance between peaks. The two components previously needed have been reduced to one.
Among the four panels in the first row, these PNS1 projections explain the horizontal variation of the original bimodal functions best, as they are almost identical to the 
raw horizontal warps of the Karcher mean,
shown in the left panel of Figure \ref{FIG:toyExample}. 
Very little variability is left for the second PNS component to explain. This suggests that the horizontal variability is almost one dimensional in some sense, which is consistent with the fact that the warping function $\gamma_i$ in this toy example can be summarized by a single parameter $a_i$. In particular, these were generated as $\gamma_i(t) = \frac{e^{a_i t}-1}{e^{a_i}-1}$, for $a_i\in [-5,5]$.

\begin{figure}[htbp]
\center
\includegraphics[keepaspectratio=true,
scale=0.5,page=1]{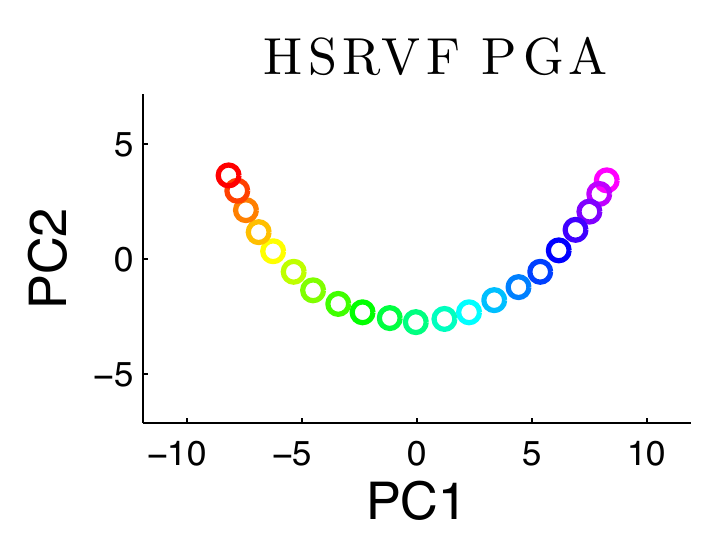}
\includegraphics[keepaspectratio=true,
scale=0.5,page=1]{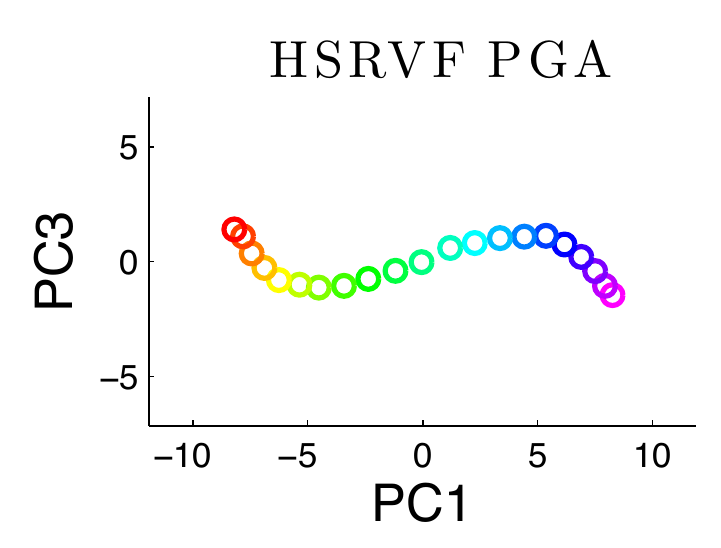}
\includegraphics[keepaspectratio=true,
scale=0.5,page=1]{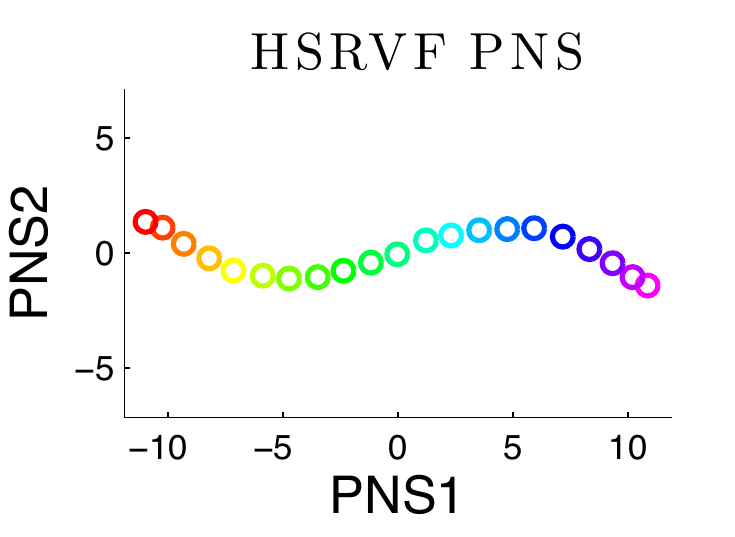}
\hskip 0.2cm
\includegraphics[keepaspectratio=true,
width=0.3\textwidth,page=1]{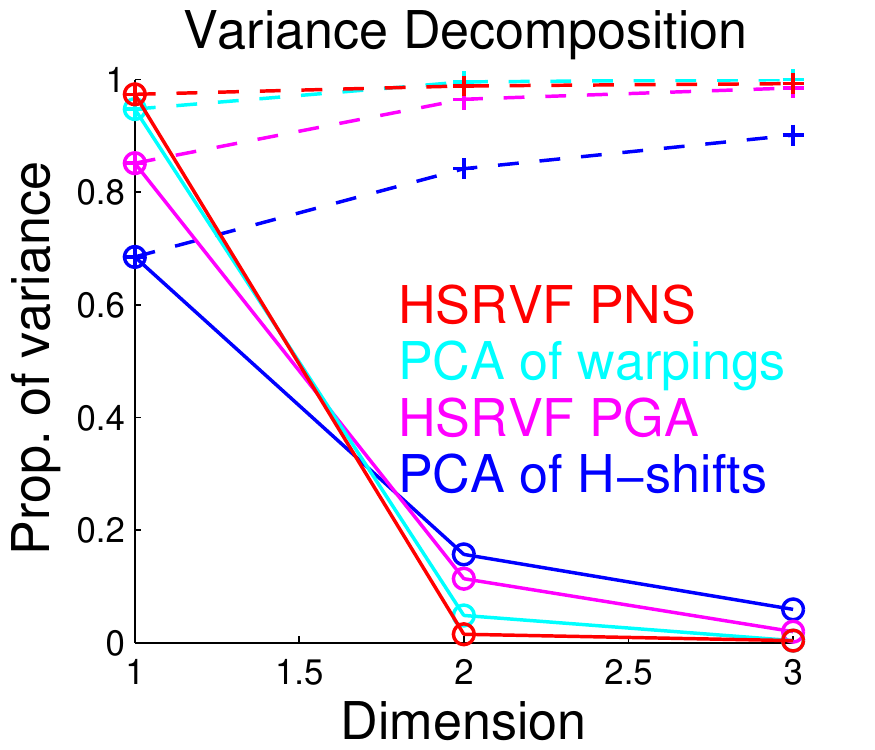}
\caption{The first two panels show scatter plots of the first three tangent PC scores in the PGA of the horizontal SRVFs. The third panel shows scatter plot of the first two PNS scores.  The aspect ratio of these plots is 1. The color is consistent with that in Figure \ref{FIG:toyExample} (left). The right panel shows scree plots of four horizontal analyses, showing both the individual (solid line) and the cumulative (dashed line)  proportion of variance explained by the first three components. The text indicates the names of the analyses, and the order of the names (from top to bottom) is consistent with the order of the proportion of variance explained by the first component (the first point). This plot shows improved signal compression by the more sophisticated methods.}
\label{FIG:scores}
\end{figure}

For further insight of this type, Figure \ref{FIG:scores} shows the score scatter plot of the first two PNS components (third panel), which has a similar pattern to that of the scatter plot of the second and the third PG scores (second panel). This is because most of the horizontal SRVF variation is along some small circle, which is captured by the PNS1. However, PGA needs two principal geodesics (the first two; see the first panel for the corresponding score scatter plot) to capture the curvature of this small circle.
In other words,  the first two PNS components are able to explain the data variability captured by the first three components in the PGA.  Thus, the PNS approach gives better signal compression by using more flexible components.

The right panel in Figure \ref{FIG:scores} visualizes the proportion of variance explained by the first three components in each of the four analyses. It is seen that the PNS of the horizontal SRVFs (red, with the highest first point) and the FPCA of the warping functions (cyan, with the second highest first point) are the top two most efficient approaches, in the sense of explaining a higher proportion of data variability using a lower number of components. However, considering the difficulty in interpreting the FPCA decomposition (see Section \ref{SEC:FPCAWarp}), PNS is the best approach for the horizontal analysis.

%================================================================

\section{Conclusions} \label{SEC:conclusion}

This paper aimed at finding an appropriate method for horizontal analysis of functional data, where the horizontal variation is separated from the vertical variation using a domain-warping method based on the Fisher Rao metric. 
Four different approaches have been discussed, including two conventional FPCA approaches (FPCA of the horizontally shifted functions of the Karcher mean and FPCA of the warping functions) and two manifold approaches based on the spherical structure of the horizontal SRVFs (PNS and PGA). A toy example of big horizontal variation is used to compare these approaches. The manifold approaches are generally better than the FPCA approaches, and the PNS works the best in terms of both the signal compression and the interpretability of the results. As a conclusion, we propose to use the PNS approach for horizontal analysis, especially when the horizontal variability is large.

%================================================================

\nocite{*}
\bibliographystyle{plain}

\bibliography{paper}

\begin{thebibliography}{10}

\bibitem{FRdist4}
N.~N. Cencov.
\newblock {\em Statistical Decision Rules and Optimal Inferences}, volume~53.
\newblock AMS, 1982.

\bibitem{PGA}
P.~T. Fletcher, C.~Lu, S.~M. Pizer, and S.~Joshi.
\newblock Principal geodesic analysis for the study of nonlinear statistics of
  shape.
\newblock {\em IEEE Trans. Medical Imaging}, 23:995--1005, 2004.

\bibitem{GG}
D.~Gervini and T.~Gasser.
\newblock Self-modeling warping functions.
\newblock {\em Journal of the Royal Statistical Society, Ser. B}, 66:959¨C971,
  2004.

\bibitem{pns}
S.~Jung, I.~L. Dryden, and J.~S. Marron.
\newblock Analysis of principal nested spheres.
\newblock {\em Biometrika}, 99(3):551--568, 2012.

\bibitem{KR}
A.~Kneip and J.~O. Ramsay.
\newblock Combining registration and fitting for functional models.
\newblock {\em Journal of American Statistical Association}, 103(483), 2008.

\bibitem{KT}
A.~Kneip and T.Gasser.
\newblock Statistical tools to analyze data representing a sample of curves.
\newblock {\em The Annals of Statistics}, 20:1266--1305, 1992.

\bibitem{FR_liu}
X.~Liu and H.~G. Mueller.
\newblock Functional convex averaging and synchronization for time-warped
  random curves.
\newblock {\em Journal of American Statistical Association}, 2004.

\bibitem{LM}
N.~Locantore, J.~S. Marron, D.~G. Simpson, N.~Tripoli, J.~T. Zhang, and K.~L.
  Cohen.
\newblock Robust principal component analysis for functional data.
\newblock {\em Test}, 8:1--73, 1999.

\bibitem{BackwardPCA}
J.~S. Marron, S.~Jung, and I.~L. Dryden.
\newblock Speculation on the generality of the backward stepwise view of pca.
\newblock Proceedings of MIR 2010: 11th ACM SIGMM International Conference on
  Multimedia Information Retrieval, Association for Computing Machinery, Inc.,
  Danvers, MA, 227-230, 2010.

\bibitem{PGAskeletal}
S.~M. Pizer, S.~Jung, D.~Goswami, X.~Zhao, R.~Chaudhuri, J.~N. Damon,
  S.~Huckemann, and J.~S. Marron.
\newblock Nested sphere statistics of skeletal models.
\newblock To appear in Proc. Dagstuhl Workshop on Innovations for Shape
  Analysis: Models and Algorithms. Springer Lecture Notes in Computer Science,
  2011.

\bibitem{RL}
J.~O. Ramsay and X.~Li.
\newblock Curve registration.
\newblock {\em Journal of the Royal Statistical Society, Ser. B}, 60:351--363,
  1998.

\bibitem{FDA}
J.~O. Ramsay and B.~W. Silverman.
\newblock {\em Functional Data Analysis}.
\newblock Springer, 2nd edition, 2005.

\bibitem{ramsay2002applied}
J.O. Ramsay and B.W. Silverman.
\newblock {\em Applied functional data analysis: methods and case studies},
  volume~77.
\newblock Springer New York:, 2002.

\bibitem{FRdist}
C.~R. Rao.
\newblock Information and accuracy attainable in the estimation of statistical
  parameters.
\newblock {\em Bulletin of Calcutta Mathematical Society}, 37:81--91, 1945.

\bibitem{FR}
A.~Srivastava, W.~Wu, S.~Kurtek, E.~Klassen, and J.~S. Marron.
\newblock Statistical analysis and modeling of elastic functions.
\newblock {\em arXiv:1103.3817}, 2011.

\bibitem{TM}
R.~Tong and H.~G. Mueller.
\newblock Pairwise curve synchronization for functional data.
\newblock {\em Biometrika}, 95(4):875--889, 2008.

\bibitem{ooda}
H.~Wang and J.~S. Marron.
\newblock Object oriented data analysis: Sets of trees.
\newblock {\em The Annals of Statistics}, 35(5):1849--1873, 2007.

\end{thebibliography}
%================================================================

\end{document}